\documentclass {amsart}

\usepackage{amsthm,amsmath,amssymb,verbatim}  
\usepackage[all]{xy}
\UseComputerModernTips
\usepackage{enumitem}

\theoremstyle{plain} 
	\newtheorem{thm}{Theorem}[section] 
	\newtheorem{prop}[thm]{Proposition} 
	\newtheorem{lemma}[thm]{Lemma}  
	\newtheorem{cor}[thm]{Corollary} 
 
\theoremstyle{definition}

	\newtheorem{examples}[thm]{Example} 
	
	\newtheorem{conjecture}[thm]{Conjecture}

\newcommand{\aro}{\longrightarrow}

\newcommand{\ep}{\varepsilon}
\newcommand{\id}{\text{id}}
\newcommand{\N}{\mathbb N}
\newcommand{\R}{\mathbb R}

\newcommand{\Lip}{\text{\rm{Lip}}}

\newcommand{\diam}{\text{\rm{diam}}}
\newcommand{\card}{\text{\rm{card}}}

\newcommand{\A}{\mathcal A}
\newcommand{\Z}{\mathbb Z}
\newcommand{\g}{\gamma}
\newcommand{\h}{\mathcal H}

\begin{document} 
 
\title	 
	{Ultrametric Spaces And the Logarithmic Ratio}
 
 
\author[H. Movahedi-Lankarani]
	{H. Movahedi-Lankarani}
\address{Department of Mathematics and Statistics\\
	Penn State Altoona\\
         Altoona, PA 16601-3760} 
\email{hxm9@psu.edu}


\keywords{logarithmic ratio, Assouad dimension, metric dimension, ultrametric, H\"older, Lipschitz} 
\subjclass{Primary 54E40; Secondary 54F45, 54F50}

\date{\today}

\maketitle 

\abstract  
 It is shown that if a compact metric space $(X, d)$ is bi-H\"older equivalent to an ultrametric space, then  the logarithmic ratio $R(X,d)$ is finite.  Conversely, if  the logarithmic ratio $R(X,d)$ is finite  and  ${\A}^*_p (X) \ne \emptyset$ for some $p \in (1, \infty )$, then $(X, d)$ is bi-H\"older equivalent to an ultrametric space.  

 It is also shown that for any $s \in [0, \infty]$ there exists a compact countable metric space $(X, d)$ with a unique cluster point such that the logarithmic ratio  $R(X, d)$ is equal $s$.  Moreover, we  prove  a bi-H\"older embedding  result for a certain class of compact totally disconnected metric spaces.
 \endabstract


\section{Introduction}

A map $f: X\to X^\prime$ of metric spaces is said to be H\"older of exponent $s>0$ if there is a constant $c>0$ such that $d^\prime (f(x_1), f(x_2))\leq c \,d(x_1,x_2)^s$ for all $x_1, x_2\in X$. If, in addition, the inverse of $f$ on $f(X)$ exists and is H\"older of exponent $t>0$, then we say that $f$ is bi-H\"older (of exponent $s$ above and $1/t$ below). 
A  bi-Lipschitz map is a bi-H\"older map of exponent 1 above and below.

 In this paper   we prove two results. The first one (Theorem 2.6) characterizes compact metric spaces which are bi-H\"older equivalent to ultrametric spaces.  We recall that an {\it ultrametric space\/} is a metric space $X$ whose metric $d$ satisfies the strong triangle inequality $d(x_1,x_2)\,\leq \max\{ d(x_1,x_3),\, d(x_2,x_3)\}$ for all $x_1,x_2,x_3\in X$.  Our second result (Theorem 4.2) deals with  bi-H\"older embeddings into Euclidean space.  Our main technical tool is the {\it{logarithmic ratio}} introduced in \cite{lograt}; see Section 2 for  definition. The logarithmic ratio is not an invariant of bi-H\"older maps, but a bi-H\"older image of a space of finite logarithmic ratio has finite logarithmic ratio. We also establish some further properties of the logarithmic ratio.  Most of our   results assume the existence of certain family of clopen partitions of $X$.  Specifically, for   each compact metric space $(X, d)$ and each real number $p\in (1, \infty)$ we define the set ${\A}^*_p (X)$ consisting of all nested sequences of clopen partitions of $X$ whose diameters do not decrease faster than an exponential in $p$ and whose logarithmic ratios converge to the logarithmic ratio of $(X, d)$. (See Section 2 for precise definition.)

The paper is organized as follows.  In Section 2 we first recall the definition and some properties of the logarithmic ratio from \cite{lograt}, and then proceed to prove  Theorem 2.6. Specifically,  we show that if a compact metric space $(X, d)$ is bi-H\"older equivalent to an ultrametric space, then logarithmic ratio $R(X,d)$ is finite.  Conversely, if the logarithmic ratio $R(X,d)$ is finite and ${\A}^*_p(X) \ne \emptyset$ for some $p \in (1, \infty)$, then $(X, d)$ is bi-H\"older equivalent to an ultrametric space. It then follows (Corollary 2.7) that a compact metric space satisfying both of these conditions is bi-H\"older embeddable as a dense subset into the inverse limit of a one-sided doubly infinite sequence of finite spaces with a comparison ultrametric. 

 In Section 3 we first show (Theorem 3.2) that for every $s\in [0,\infty]$ there exists a compact countable metric space with a unique cluster point  whose logarithmic ratio is equal to $s$.  We next  show that  the logarithmic ratio  respects inclusion (Proposition 3.3)  and behaves reasonably well with respect to  Cartesian products.

Section 4 is devoted to finite dimensional bi-H\"older embeddings. The aim is to establish embedding results similar to the one in \cite{LU-ML} for ultrametric spaces for other classes of totally disconnected spaces. The problem of characterizing all metric spaces which admit bi-H\"older embeddings  into $\Bbb R^N$ for some $N\in \Bbb N$ is still open. However, there has been notable progress. Indeed, in \cite{assouad2, assouad3} Assouad defines a bi-Lipschitz  invariant called the {\it{Assouad dimension}}, and denoted by $\dim_A$, (see Section 4) which is a generalization of the dimensional order of Bouligand \cite{b}. ouligandBy using this notion of dimension, in combination with his work in \cite{assouad1}, he proves that if $(X,d)$ is a metric space of finite Assouad dimension and $s\in (0,1)$, then there is a bi-Lipschitz embedding of $(X,d^s)$ into $\Bbb R^N$ for {\it some\/} $N\in\Bbb N$.
In Theorem 4.2  we give sufficient conditions for existence of a bi-H\"older embedding $f$ into $\R^N$, where both $f$ and $f^{-1}$ are H\"older of {\it same} exponent $s>0$.  Furthermore, we   give  an estimate on the size of the integer $N$.  Our  result is for a certain class of compact metric spaces of finite (but $>1$) logarithmic ratio $R$ and finite Assouad dimension $D$. The estimate for $N$ depends on both $R$ and $D$ which in a special case reduces to $N > (D + R -1)(2 R + 1)$. 
We conclude the paper with a few remarks in Section 5.  The work in this paper is essentially a continuation of \cite{LU-ML} and \cite{lograt}.

Finally, since we deal only with compact metric spaces, we assume that all our spaces have diameter $\le 1$.  Also, the author thanks Jouni Luukkaninen and an anonymous reader for their careful reading  of  earlier versions of this paper and making several corrections and suggestions.



\section{Spaces bi-H\"older equivalent to ultrametric spaces}

We begin this section by recalling the definition and some properties of the logarithmic ratio from \cite{lograt}.  To this end, let $\alpha = \{ A_1, A_2,..., A_n \}$ be a {\it clopen\/} partition of a compact metric space $(X,d)$ with $n = \card\, \alpha$. Define the {\it diameter\/} $\delta (\alpha) \in [0, {\infty})$ of $\alpha$ by setting $\delta (\alpha) =  \max \{\diam    (A_i) \bigm\vert 1\leq i\leq n\}$ if $n\geq 1$ and $\delta (\alpha ) = 0$ if $n = 0$. Define the {\it gap\/} $\gamma (\alpha) \in [0, \infty )$ of $\alpha$ by setting $\gamma (\alpha) =  \min \{ d(x,y) \bigm\vert  x \in A_i,\, y \in A_j,\, i\neq j \}$ if $n > 1$ and by $\gamma (\alpha ) = \diam (X)$ if $n\leq 1$ (with the convention that $\diam (\emptyset)  = 0$). Furthermore, define the {\it logarithmic ratio} $R (\alpha) \in [0,\infty ]$ of $\alpha$ by setting $R (\alpha) = \log \gamma (\alpha) /\log \delta (\alpha) $ if $0 < \delta(\alpha) < 1$ and $(0 <)\,\gamma (\alpha) < 1$; by $R(\alpha) = 0$ if $\delta (\alpha) = 0$; and by $R( \alpha) = +\infty$ if either $0 < \delta (\alpha) < 1$ and $\gamma (\alpha) \geq 1$ or $\delta (\alpha) \geq 1$. Then the {\it logarithmic ratio} of $X$ (with respect to the metric $d$) is defined by

\[R(X,d) = \lim_{r\to 0} \, \inf \left\{ R(\alpha) \bigm\vert \alpha \,\text{ 
is a clopen partition of}\,\, X \,\text {and}\, \,\delta (\alpha) <r \right\}\]

\noindent with the convention that  $ \inf \emptyset = +\infty$. The following results are proved in \cite{lograt}:

\begin{enumerate}
\item If $X$ is a finite set, then $R(X,d) = 0.$
\item If $R(X,d)<\infty$, then $(X,d)$ is totally disconnected.
\item If $(X,d)$ and $(X^\prime,d^\prime)$ are bi-Lipschitz  equivalent, then $R(X,d)= R(X^\prime,d^\prime).$
\item If $(X,d)$ is a compact ultrametric space, then $R(X,d) \leq 1$. 
\item Let $R(X,d)\leq 1$ and suppose that there is a sequence $\bigl(\alpha_n\bigr)_{n\geq 1}$ of clopen partitions of $X$ with

(a) $\alpha_{n+1}$ is a refinement of $\alpha_n$ for all $n\geq 1$, and $0 < \delta(\alpha_n) \rightarrow 0$ as $n \to \infty$,

(b) There is a constant $C\in [1,\infty)$ such that $\delta (\alpha_n)/\delta (\alpha_{n+1})\leq C$ for all $n\geq 1$,  and

(c) $\lim_{n\to \infty} R(\alpha_n) = R(X,d)$.

\noindent Then $(X,d)$ is bi-H\"older equivalent to an ultrametric space.
\end{enumerate}

\noindent Two points $x_1, x_2 \in X$ are said to be {\it associated endpoints} of $X$ if there is a clopen subset $A \subset X$ such that $x_1 \in A$, $x_2 \in X \setminus A$, and for the partition $\alpha = \{A, X \setminus A\}$ we have $\gamma (\alpha) = d (x_1, x_2)$. Then the {\it largest gap}  in $X$ is defined to be $\Gamma (X) = \sup\{d(x, x^{\prime}) \bigm\vert x, x^{\prime} \in X  {\text{ are associated endpoints of }} X\}$.

\begin{enumerate}[resume]
\item   Suppose that there exists a sequence $\bigl(\alpha_n\bigr)_{n\geq 1}$ of clopen partitions of $X$ with

(a)  $0 < \delta(\alpha_n) \searrow 0$ as $n \to \infty$, and

(b) There is a constant $C\in [1,\infty)$ such that for all $n \ge 1$ there is $A \in \alpha_n$ with $\diam (A) = \delta (\alpha_n)$ and  $\Gamma (A) \le C \gamma (\alpha_{n +1 })$.
  
\noindent Then $R(X,d) = \liminf_{n\to \infty} R(\alpha_n)$.
\end{enumerate}

\noindent With minor modifications, a proof similar to that of (3) above yields the following result.

\begin{lemma}
 Let $(X,d)$ and $(X^\prime,d^\prime)$ be compact metric spaces. If there is a bi-H\"older equivalence $(X,d) \aro (X^\prime,d^\prime)$ of exponent {\it s} above and {\it t} below, $0<s \leq t$, then 
\[\frac {s}{t}\,R(X,d) \,\leq \, R(X^\prime,d^\prime) \,\leq \,\frac{t}{s}\, R(X,d).\]
\end{lemma}

\begin{cor} Let $(X, d)$ be a compact metric space and let $0 < s \le 1$. Then $R(X,d^s) = R(X, d)$.
\end{cor}

\noindent  By combining Lemma 2.1 with statement (4) above we see that if a compact metric space $(X,d)$ is bi-H\"older equivalent to an ultrametric space, then $R(X,d) < \infty$. There is also a partial converse to this statement. Specifically, let $(X,d)$ be a compact metric space and set, cf. \cite{ultgeomeas},
\[ {\A} (X) = \left\{\bigl(\alpha_n\bigr)_{n\geq 1} \bigm\vert \alpha_n \text{ is a clopen partition of } X \text{ and } \lim_{n \to \infty} \delta (\alpha_n) = 0 \right\}.\]

\noindent Clearly, ${\A}(X) \ne \emptyset$ if and only if $X$ is totally disconnected. Also, by \cite[Lemma 2.2]{lograt}, if $(X,d)$ is nondiscrete, then $\lim_{n\to\infty}\gamma (\alpha_n) = 0$ for every sequence $\bigl(\alpha_n\bigr)_{n\geq 1} \in {\A} (X)$. Clearly then, if $(X, d)$ is totally disconnected and nondiscrete, exactly one of the following two properties hold:

\begin{enumerate}

\item  There is a sequence $\bigl(\alpha_n\bigr)_{n\geq 1} \in {\A} (X)$ with $\gamma (\alpha_n) > \delta (\alpha_n)$ for all $n \in \N$.

\item Each sequence $\bigl(\alpha_n\bigr)_{n\geq 1} \in {\A} (X)$ has a subsequence $\bigl(\alpha_{n_k}\bigr)_{k\geq 1} \in {\A} (X)$ with $\gamma (\alpha_{n_k}) \le \delta (\alpha_{n_k})$ for all $k \in \N$.
\end{enumerate}

\noindent Let 
\[\hat{\A} (X) = \left\{\bigl(\alpha_n\bigr)_{n\geq 1} \in {\A} (X) \bigm\vert  \alpha_{n+1} {\text{is a refinement of }} \alpha_n {\text{ and}} \liminf_{n\to\infty} R(\alpha_n) = R(X,d)\right\}.\]

\noindent  If $\alpha$ and $\beta$ are clopen partitions of $X$ with $\delta (\beta) < \gamma (\alpha)$, then $\beta$ is a refinement of $\alpha$.   Hence, if ${\A}(X) \ne \emptyset$, then $\hat{\A} (X)\ne \emptyset$. Consequently, if $R(X,d) < \infty$, then $\hat{\A} (X)\ne \emptyset$.  Furthermore, if $1 < R(X, d) < \infty$, then there is a sequence $\bigl(\alpha_n\bigr)_{n\geq 1}\in \hat{\A} (X)$ with $\g(\alpha_n) < \d(\alpha_n)$. Similarly, if $R(X, d) < 1$, then there is a sequence $\bigl(\alpha_n\bigr)_{n\geq 1}\in \hat{\A} (X)$ with $\d(\alpha_n) < \g(\alpha_n)$. For $p\in (1, \infty)$ we set

\[\hat{\A}_p (X) = \left\{\bigl(\alpha_n\bigr)_{n\geq 1}\in \hat{\A}  (X) \bigm\vert   {\text{there is }} a > 0 {\text{ with }} a \, \delta (\alpha_n)^p\,\leq\,\delta (\alpha_{n+1})\,\leq\,\delta (\alpha_n)  \text{ for all } n\geq 1  \right\}\]

\noindent and we note that $\hat{\A}_p (X) \subset \hat{\A}_q (X)$ for $q\geq p$.  However, the nested family of sets

\[\emptyset\subset\dots\subset \hat{\A}_p (X)\subset\dots\subset \hat{\A}_q (X)\subset\dots\subset \hat{\A} (X),\  \  1 < p \leq q\]

\noindent   is not necessarily a filtration because $\bigcup_{1 < p} \hat{\A}_p (X) \ne \hat{\A} (X)$ is possible. In general, for a given $p\in (1, \infty)$, the set $\hat{\A}_p (X)$ may be empty.  Indeed, the condition

\[\hat{\A}_p(X)\ne\emptyset \ {\text { for some }}\   p\in (1, \infty)\]

\noindent is invariant under a certain  family  of bi-H\"older maps;  see also examples 3.1 and 6.1 below.   

\begin{lemma} 
Let $(X, d)$ be  a compact metric space.   Then for each  $p > 1$, the condition      $\hat{\A}_p(X) \ne \emptyset$  is invariant under bi-H\"older maps of the  same  exponent above as below; in particular,  bi-Lipschitz maps.
\end{lemma}

\begin{proof} 
 Let  $f \: (X,d) \aro (X^\prime,d^\prime)$  be a bi-H\"older map of exponent {\it s} above and {\it t} below with $0<s \leq t$.  Hence, we have 

\[c_1 \, d (x_1, x_2)^t \le d^{\prime} \bigl(f (x_1), f (x_2)\bigr) \le c_2\,  d (x_1, x_2)^s\]
for some constants  $c_1 > 0$ and $c_2 > 0$ and all $x_1, x_2 \in X$.  Assume that $\hat{\A}_p (X) \ne \emptyset$ for some $p \in (1, \infty)$.   Then there is $\bigl(\alpha_n\bigr)_{n\geq 1}\in \hat{\A}  (X)$ and $a > 0$ with  $a \, \delta (\alpha_n)^p\,\leq\,\delta (\alpha_{n+1})\,\leq\,\delta (\alpha_n)$   for all  $n\geq 1$. Also,   for each $n \ge 1$,   $f (\alpha_n) = \bigl\{f (A) \bigm\vert A \in \alpha_n\bigr\}$ is a clopen partition of $X^{\prime}$ with 
$c_1 \gamma (\alpha_n)^t \le \gamma \bigl(f (\alpha_n)\bigr) \le c_2 \gamma (\alpha_n)^s$ and  $c_1 \delta (\alpha_n)^t \le \delta \bigl(f (\alpha_n)\bigr) \le c_2 \delta (\alpha_n)^s \to 0$ as  $n \to \infty$.  Moreover, for each $n \ge 1$, $f (\alpha_{n +1})$ is a refinement of $f (\alpha_n)$ and $a^{\prime} \delta \bigl(f (\alpha_n)\bigr)^{p^{\prime}} \le \delta \bigl(f (\alpha_n)\bigr)$ with  $p^{\prime} = (t / s) \, p$  and    $a^{\prime} = (c_1 / c_2^{p^{\prime}}) \, a^t$.   If $(X, d)$ is discrete, then for all  sufficiently large $n$ we have $\delta (\alpha_n)  = 0$ and $\delta \bigl(f (\alpha_n)\bigr) = 0$ so that $R (\alpha_n) = 0 = R (X, d)$ and $R \bigl(f (\alpha_n)\bigr) = 0 = R (X^{\prime}, d^{\prime})$.  If $(X, d)$ is not discrete, we have 

\[\frac{s}{t}\, R (X, d) \le R (X^{\prime}, d^{\prime}) \le \liminf_{n \to \infty} R \bigl(f (\alpha_n)\bigr) \le \frac{t}{s}\, R (X, d). \]

\noindent This proves the  lemma.
\end{proof}

\begin{lemma}
 Let $(X, d)$ be  a compact metric space with $0 < R(X, d)  \le 1$. If $(X, d)$ is bi-Lipschitz equivalent to an ultrametric space, then $\hat{\A}_p(X) \ne \emptyset$ for any $p\in (1/R (X, d), \infty)$.
\end{lemma} 

\begin{proof} If $(X, d)$ is a compact ultrametric space with $0 < R (X, d)$, then   the metric $d$ takes values in a set $\{0\} \cup\{r_1 > r_2 > \cdots\}$ with $r_n \searrow 0$.  We recall that any two balls in an ultrametric space are either disjoint or one is contained in the other.  Hence, letting $\alpha_n$ be the clopen partition of $X$ with closed $r_n$-balls, we see that $\delta (\alpha_n) = r_n$ and $\gamma (\alpha_n) = r_{n - 1}$ so that $R (\alpha_n) = \log r_{n - 1}/\log r_n$ for $n > 1$. Also,  for $n > 1$ there is  $A  \in \alpha_n$  with $\diam (A) = r_n$  and   $\Gamma(A) = r_n  = \g(\alpha_{n + 1})$. It then follows   that

\[ R (X, d) = \lim_{k \to \infty} \, \inf_{k \le n} R (\alpha_n) = \lim_{k \to \infty} \, \inf_{k \le n} \frac{\log r_{n - 1}}{\log r_n};\]

\noindent the  lemma now follows.
\end{proof}

The case of $1 < R (X, d)$ remains open; see Conjecture 2.8 below.     Now,  let

\[{\A}^* (X) = \left\{\bigl(\alpha_n\bigr)_{n\geq 1} \in \hat{\A} (X) \bigm\vert    \liminf_{n\to\infty} R(\alpha_n) =  \lim_{n\to\infty} R(\alpha_n) = R(X,d)\right\},\] 

\noindent  and,  for $p \in (1, \infty)$, let  ${\A}^*_p  (X) = {\A}^* (X) \cap \hat{\A}_p (X)$. 
That is, if $\bigl(\alpha_n\bigr)_{n\geq 1} \in {\A}^*_p (X)$, then, for all  $n\geq 1$,   $\alpha_{n + 1}$ is a refinement of $\alpha_n$,  there is  $a > 0$  with  $a \, \delta (\alpha_n)^p\,\leq\,\delta (\alpha_{n+1})\,\leq\,\delta (\alpha_n)$,   and $\lim_{n\to\infty} R(\alpha_n) = R(X,d)$.  Moreover, we have

\[\frac{1}{p}  \le \liminf_{n \to \infty} \frac{\log \gamma_n}{\log \gamma_{n + 1}}  \le  \limsup_{n \to \infty} \frac{\log \gamma_n}{\log \gamma_{n + 1}} \le 1.\]

\begin{lemma}
 Let $(X,d)$ be a compact metric space with $R(X,d) < \infty$ and suppose that ${\A}^*_p(X)\ne\emptyset$ for some $p\in (1,\infty)$. Then there is a compatible ultrametric $\rho$ on $X$ such that $(X,\rho)$ is bi-H\"older equivalent to $(X,d)$.
\end{lemma}

\begin{proof}
If $X$ is discrete, the lemma is trivially true. Otherwise,  let $R(X,d) = R < \infty$ and let $\bigl(\alpha_n\bigr)_{n\geq 1}\in {\A}^*_p(X)$ be a sequence of partitions of $X$. Then  we have

\begin{enumerate}
\item $\alpha_{n+1}$ is a refinement of $\alpha_n$ for all $n\geq 1$;
\item $\delta(\alpha_n) \to 0$ and $\gamma(\alpha_n) \to 0$;
\item $\lim_{n\to\infty}R(\alpha_n)\,=\,R(X,d)\,=\,R$;
\item $ a \, \delta(\alpha_n)^p \,\leq\, \delta(\alpha_{n+1})$ for some $a > 0$ and all $n\geq 1$.
\end{enumerate}

\noindent Furthermore, given $\ep >0$, property (3) implies that there is some fixed $m \geq 1$ such that for $n > m$  we have

\[\delta(\alpha_n)^{R+\varepsilon} < \gamma(\alpha_n) < \delta(\alpha_n)^{R-\varepsilon}. \] 

Now let $\alpha_0 = \{X\}$. Let $x_1$ and $x_2$ be distinct points in $X$ and set $\rho (x_1,x_2)\,=\,\delta(\alpha_n)$ if $x_1$ and $x_2$ are in the same member of $\alpha_n$ but in different members of $\alpha_{n+1}$. Then it is clear that $\rho$ is an ultrametric on $X$ and that it is compatible. Also, we have
\[d(x_1,x_2) \leq \delta(\alpha_n) \,=\,\rho (x_1,x_2).\]

\noindent To get an inequality on the lower side, let $x_1$, $x_2$, and $n \ge m$ be as above and observe that

\[d (x_1,x_2) \geq\gamma(\alpha_{n+1}) > \delta(\alpha_{n+1})^{R+\varepsilon} \geq a^{R + \varepsilon} \delta(\alpha_n)^{p(R+\varepsilon)}\,=\, a^{R + \varepsilon} \rho (x_1,x_2)^{p(R+\varepsilon)}.\]

\noindent Hence, if $x_1$ and $x_2$ are in the same member of $\alpha_m$, we have

\[a^{R + \varepsilon} \rho (x_1,x_2)^{p (R+\varepsilon)} \leq d(x_1,x_2) \leq \rho (x_1,x_2).\]

\noindent If $x_1$ and $x_2$ are not in the same member of $\alpha_m$, then

\[\rho (x_1, x_2)^{p (R + \ep)} \le \delta (\alpha_0)^{p (R + \ep)} \le \, \frac{\delta (\alpha_0)^{p (R + \ep)}}{\gamma (\alpha_m)} \, d(x_1, x_2).\]

\noindent Thus, if $K = \min \{a^{R + \varepsilon},  \gamma (\alpha_m) \delta (\alpha_0)^{- p (R + \ep)}\}$, then

\[K \rho (x_1, x_2)^{p (R + \ep)} \le d (x_1, x_2) \le \rho (x_1, x_2)\]

\noindent for all $x_1, x_2 \in X$.
\end{proof}

In the above lemma if $ R(X,d)\,=\,R$, then $R /p (R+\varepsilon) \leq R(X,\rho) \leq 1$. Moreover, we have\footnote{Throughout, $\dim_H$ denotes the Hausdorff dimension.} $\dim_H(X,d)\,\leq\, \dim_H(X,\rho)\,\leq\, p(R+\varepsilon) \dim_H(X,d)$.    We note that the above lemma is  related to  \cite[Proposition 2.12]{lograt}. As an immediate corollary of this lemma along with Lemma 2.1 we get the following result.

\begin{thm} 
\label{ult}
Let $(X,d)$ be a compact metric space. 
\begin{enumerate}
\item If $(X,d)$ is bi-H\"older equivalent to an ultrametric space, then $R(X,d) < \infty$.
\item If $R(X,d) < \infty$ and ${\A}^*_p(X)\ne\emptyset$ for some $p\in (1, \infty)$, then $(X,d)$ is bi-H\"older equivalent to an ultrametric space.
\end{enumerate} 
\end{thm}

It is shown in \cite[Theorem 2.2 and Remark 2.3.1]{LU-ML} that every compact ultrametric space is bi-Lipschitz embeddable as a dense subset into the inverse limit of a   sequence of  finite discrete spaces with a so-called comparison ultrametric. This result in combination with Theorem \ref{ult} yields the following corollary.

\begin{cor} 
Let $(X,d)$ be a compact metric space with $R(X, d) < \infty$ and suppose that ${\A}^*_p(X)\ne\emptyset$ for some $p\in (1, \infty)$. Then $(X,d)$ is bi-H\"older embeddable as a dense subset into the inverse limit of a   sequence of  finite  discrete spaces with a comparison ultrametric.
\end{cor}

It follows that if Conjecture 2.8 below is true, then the class of compact ultrametric spaces constructed in \cite[Construction 2.1]{LU-ML} is in a sense universal for the category of compact metric spaces of finite logarithmic ratio and bi-H\"older embeddings.

\begin{conjecture} 
If $(X,d)$ is a compact metric space with $0 < R(X,d) < \infty$, then ${\A}^{\ast}_p (X) \ne \emptyset$ for some $p\in (1,\infty )$.
\end{conjecture}



\section{More on the logarithmic ratio}

In this section we establish some further properties  of the logarithmic ratio.  We begin with the following class of examples.

\begin{examples}  Let $(r_n)_{n \geq 1}\subset (0,1]$ be a strictly decreasing sequence with limit 0 and with  $r_{n - 1} -  r_n \le r_n + r_{n + 1}$ for all $n > 1$. Set $X\, =\, \{0\}\cup \{ r_n \mid n\geq 1\}$ with the inherited metric $\mid \cdot \mid$. For each $n > 1$, let $\alpha_n$ be the partition of $X$ given by

\[\alpha_n \, =\, \left\{ X \setminus \{r_1, r_2, \cdots, r_{n - 1}\}, \{r_1\}, \{r_2\}, \dots ,\{r_{n-1}\} \right\}.\]

\noindent Then, for each $n>1$, we have $\delta (\alpha_n) = r_n$ and $\gamma (\alpha_n) = r_{n - 1} - r_n$ so that $R(\alpha_n)\,=\, \log (r_{n-1}\,-\,r_n)/\log r_n$. Again,  for $n>1$ and $A = X \setminus \{r_1, r_2, \cdots, r_{n - 1}\} \in \alpha_n$,  we have $\diam (A) = r_n = \delta (\alpha_n)$ and  the largest gap in $A$ is $\Gamma (A) = r_n - r_{n + 1} = \gamma (\alpha_{n + 1})$.  It now follows   that

\[R (X, \mid \cdot \mid) = \lim_{k \to \infty} \, \inf_{k \le n} R (\alpha_n) = \lim_{k \to \infty} \, \inf_{k \le n} \frac{\log(r_{n-1}- r_n)}{\log r_n}.\]

\noindent We consider five cases:

\vskip5pt

\noindent {\it{Case 3.1.1.}} Let $r_n = 2^{- n!}$. Then $R (X, \mid \cdot \mid) = 0$. This is an example of a {\it nondiscrete} compact metric space with {\it zero} logarithmic ratio.\footnote{See also Example 5.1.} Furthermore, for each $p > 1$ and each $n(p) \in \N$  we have $\lim_{n \to \infty} R (\alpha_n) = 0$ but  $\bigl(\alpha_n\bigr)_{n\geq n(p)} \notin {\A}^*_p(X)$.\footnote{This, of course,  is far from saying  that ${\A}^*_p (X) = \emptyset$ for all $p > 1$.} We note, however, that $(X, \mid \cdot \mid)$ is  bi-Lipschitz equivalent to an ultrametric space. Indeed, let $\rho$ be the ultrametric on $X$ defined by setting

\[\rho (r_n, r_n) = 0, \ \ \  \rho (0, r_n) = r_n, \ \ {\text{ and }}\ \rho (r_m, r_n) = \max \{ r_m, r_n\} {\text{ for }} m \ne n.\]

\noindent Then we have $\mid \cdot \mid \, \le \rho \le 2 \mid \cdot \mid$.

\vskip5pt

\noindent {\it{Case 3.1.2.}} For $0 < s < 1$, let $r =  2^{- s^2 / (1 - s)}$ and let $r_n = r^{1/s^n}$. Then we have $\lim_{n \to \infty} R (\alpha_n) = s$ so that $R (X, \mid \cdot \mid) = s$. Furthermore,  $(X, \mid \cdot \mid)$ is  bi-Lipschitz equivalent to an ultrametric space. Indeed, if we let $\rho$ be the ultrametric on $X$ defined in Case 3.1.1, then we have $\mid \cdot \mid \, \le \rho \le \left(2^{1/s}/(2^{1/s} - 1)\right) \mid \cdot \mid$.   If $p > 1/s$, then for  $n \ge 1$ we have $\delta (\alpha_n)^p \le \delta (\alpha_{n + 1})$ so that $\bigl(\alpha_n\bigr)_{n \ge 1}  \in {\A}^*_p(X)$. 

\vskip5pt

\noindent {\it{Case 3.1.3.}} Let $r_n = 2^{-n}$. Then $\lim_{n \to \infty} R (\alpha_n) = 1$ so that $R (X, \mid \cdot \mid) = 1$. Furthermore, if $p > 1$, then $2^{p - 2} \delta (\alpha_n)^p \le \delta (\alpha_{n + 1})$ so that $\bigl(\alpha_n\bigr)_{n \ge 1}  \in {\A}^*_p(X)$ implying that $(X, \mid \cdot \mid)$ is  bi-H\"older equivalent to an ultrametric space. Indeed, $(X, \mid \cdot \mid)$ is  bi--Lipschitz equivalent to an ultrametric space, for if we let $\rho$ be the ultrametric on $X$ defined in Case 3.1.1, then we have  $\mid \cdot \mid \, \le \rho \le 2 \mid \cdot \mid$.

\vskip5pt

\noindent {\it{Case 3.1.4.}} For $1 < s < \infty$, set $r_n = n^{1/(1-s)}$. Then we  have $\lim_{n \to \infty} R (\alpha_n) = s$ so that $R (X, \mid \cdot \mid) = s$. Furthermore, if $p > 1$, then $2^{1/(1 - s)} \delta (\alpha_n)^p \le \delta (\alpha_{n + 1})$ so that $\bigl(\alpha_n\bigr)_{n \ge 1}  \in {\A}^*_p(X)$  implying that $(X, \mid \cdot \mid)$ is  bi-H\"older equivalent to an ultrametric space.

\vskip5pt

\noindent {\it{Case 3.1.5.}} For $n \ge 3$,  let $r_n = 1/\log n$.     Then  we have $\lim_{n \to \infty} R (\alpha_n) = \infty$ so that $R (X, \mid \cdot \mid) = \infty$ implying that $(X, \mid \cdot \mid)$ is {\it{not}} bi-H\"older equivalent to an ultrametric space.  However, for $p > 1$, we have  $(\log 3/\log 4) \delta (\alpha_n)^p \le \delta (\alpha_{n + 1})$ so that $\bigl(\alpha_n\bigr)_{n \ge 1}  \in {\A}^*_p(X)$.  This is an example of a {\it totally disconnected} compact metric space with {\it infinite} logarithmic ratio. 
\end{examples}

\begin{thm} For any $s\in [0, \infty]$ there exists a compact countable metric space $(X,d)$ with a unique cluster point such that $R(X,d) = s$. Indeed, the space $X$ may be taken to be a strictly decreasing sequence in the interval $(0, 1]$ together with its limit point $0$.
\end{thm}

The following result is due to Jouni Luukkainen\footnote{Private Communication}.

\begin{prop} {\rm{[J. Luukkainen]}} Let $Y$ be a closed subset of a compact metric space $(X, d)$. Then $R (Y, d) \le R (X,d)$.
\end{prop}

\begin{proof} We may assume that $R (X, d) < \infty$ and that $R (Y, d) > 0$. It then follows that  both $X$ and $Y$ are totally disconnected and nondiscrete. Let $\bigl(\alpha_n\bigr)_{n \geq 1} \in {\A}^*(X)$ with $\delta (\alpha_n) \in (0, 1)$ and $\gamma (\alpha_n) \in (0, 1)$ for all $n \geq 1$. Then $\delta (\alpha_n) \to 0$, $\gamma (\alpha_n) \to 0$, and $R (\alpha_n) = \log \gamma (\alpha_n)/ \log \delta (\alpha_n) \to R (X, d)$. For each $n \geq 1$, we have that $\beta _n = \{ A\cap Y \mid A \in \alpha_n\} \setminus \{\emptyset\}$ is a clopen partition of $Y$ with $0 < \delta (\beta_n) \le \delta (\alpha_n)$ and $0 < \gamma (\alpha_n) \le \gamma (\beta_n)$. Therefore, $\delta (\beta_n) \to 0$ and so $\gamma (\beta_n) \to 0$ by \cite[Lemma 2.2]{lograt}. Hence, we may assume that $\gamma (\beta_n) < 1$ for each $n$.  Thus, for each $n \geq 1$ we have $R (\beta_n) \le R (\alpha_n)$ implying that

\[R (Y , d) \le \lim_{k \to \infty} \, \inf_{k \le n} R (\beta_n) \le \lim_{n\to\infty} R (\alpha_n) = R (X , d) . \]
 \end{proof}

Next,  given a nonempty family of metric spaces, there are many choices for the metric on the set theoretic product of the family.  Following  \cite{invlim, HGlim} we equip the set theoretic product with the $\sup$-metric $d_{\sup}$, and call it the {\it{metric product}} of the family.  Specifically, let $S \ne \emptyset$ be a set and for each $s \in S$, let $\bigl(X_s, d_s\bigr)$ be a  metric space. Then the metric product of this family is $\bigl(\prod_{s \in S}X_s, d_{\sup}\bigr)$, where

\[d_{\sup} (x, x^{\prime}) = \sup_{s \in S} d_s (x_s, x^{\prime}_s),\]

\noindent with the canonical projections $\pi_t : \prod_{s \in S} X_s \aro X_t$. We note two issues here. One: In general, $d_{\sup}$ is an extended metric because $\infty$ is allowed as a value; consider $\R^{\N}$. However, this is not a problem because we are concerned with compact metric spaces and clopen partitions with arbitrarily small diameters. It is clear that if, for instance, there is a uniform upper bound on the diameters of the family $(X_s)_{s \in S}$, then $d_{\sup}$ is a metric. And two: The topology induced by $d_{\sup}$ on an infinite product of metric spaces is not necessarily the product topology; $\bigl(\Z/2 \Z, d_{\sup}\bigr)$ is discrete.  However, the metric product satisfies the following  universal properties with respect to Lipschitz  maps:

\begin{enumerate}
\item For any $t \in S$, $g : X_t \aro Y$, and $f : Y \aro \prod_{s \in S} X_s \ne \emptyset$ we have $\Lip \,  (\pi_t) = 1$, $\Lip \,  (g \circ  \pi_t) = \Lip \,  (g)$, $\diam \bigl((g \circ \pi_t) ( \prod_{s \in S} X_s)\bigr) = \diam \bigl(g (X_t)\bigr)$, and $\diam \bigl(f (Y)\bigr) = \sup_{s \in S} \diam \bigl((\pi_s \circ f) (Y)\bigr)$.

\item Given a family of maps $\bigl(f_s : Y \aro X_s\bigr)_{s \in S}$, there is a unique map $f : Y \aro \prod_{s \in S} X_s$ with $\pi_s \circ f = f_s$ for all $s \in S$, and satisfying $\Lip \, (f) \le \sup_{s \in S} \Lip\,  (f_s)$.
\end{enumerate}

The following result is an immediate corollary of the Theorem 3.3.

\begin{cor} Let $S \ne \emptyset$ be a set and for each $s \in S$, let $\bigl(X_s, d_s\bigr)$ be a compact and nonempty metric space.   If $\bigl(\prod_{s \in S} X_s, d_{\sup}\bigr)$ is compact, then $R \bigl(\prod_{s \in S} X_s, d_{\sup}\bigr) \geq \sup_{s \in S}  R(X_s, d_s)$.
\end{cor}

The following example shows that, in general, the inequality in this corollary cannot be improved.

\begin{examples}  Let $\bigl(r_n\bigr)_{n \ge 1} \subset (0, 1]$ be a  strictly monotone decreasing sequence converging to $0$. For each $n \ge 1$, let $X_n = r_n (\Z/2\Z)$ with the obvious ultrametric $d_n$ . Then $R (X_n, d_n) = 0$ for all $n \ge 1$, and $R \bigl(\prod_{n \ge 1} r_n (\Z/2\Z), d_{\sup}\bigr) \le 1$ because the metric product of a family of ultrametric spaces  is an ultrametric space. For each $n \ge 1$, let $\alpha_n$ be the partition of $(\prod_{k \ge 1}r_k (\Z/2\Z), d_{\sup}\bigr)$ by the balls of radius $r_n$.  It follows  that $R (\alpha_n) = \log r_{n - 1}/\log r_{n}.$  Again,  for $n \ge 1$ and $A$ any ball of radius $r_n \in \alpha_n$,  the largest gap in $A$ is $\Gamma(A) = r_{n}  = \gamma (\alpha_{n + 1})$.  It then follows   that

\[R \bigl(\prod_{n \ge 1} r_n (\Z/2\Z), d_{\sup}\bigr)  = \lim_{k \to \infty} \, \inf_{k \le n} R (\alpha_n) = \lim_{k \to \infty} \, \inf_{k \le n} \frac{\log r_n}{\log r_{n + 1}} \le 1.\]

\noindent By assuming $r_1 < 1$ and setting $r_{n + 1} = r_n^{1/t}$ for any $t \in (0, 1)$ and all $n \ge 1$, we obtain 

\[R \bigl(\prod_{n \ge 1}r_n (\Z/2\Z), d_{\sup}\bigr) = t \in (0, 1). \  \square \] 
\end{examples}

What is also  true is the following:  We recall that for two subset $A$ and $A^{\prime}$ of a metric space, we set 

\[{\text{gap}} (A, A^{\prime}) = \inf \bigl\{d \, (x, x^{\prime}) \bigm\vert x \in A {\text{ and }} x^{\prime} \in A^{\prime}\bigr\}.\]

\noindent If $\alpha \ne \emptyset$ and $\alpha^{\prime} \ne \emptyset$ are clopen partitions of compact metric spaces $(X, d)$ and $(X^\prime, d^\prime)$, respectively, then $\alpha \times \alpha^{\prime} = \bigl\{ A \times A^{\prime} \mid A \in \alpha {\text{ and }} A^{\prime} \in \alpha^{\prime}\bigr\}$ is a clopen partition of $(X \times X^\prime, d_{\sup})$ with $\delta (\alpha \times \alpha^{\prime}) = \max \{\delta (\alpha), \delta (\alpha^{\prime})\}$. Moreover,  if $\card (\alpha) > 1$ and $\card (\alpha^{\prime}) >1$, then

\[\aligned
\gamma (\alpha \times \alpha^{\prime}) &= \min \bigl\{{\text{gap}} (A_1 \times A^{\prime}_1, A_2 \times A^{\prime}_2) \bigm\vert A_1, A_2 \in \alpha, A^{\prime}_1, A^{\prime}_2 \in \alpha^{\prime}, A_1 \times A_1^{\prime} \ne A_2 \times A_2^{\prime}\bigr\}\\
&= \min \left\{\max \bigl\{{\text{gap}} (A_1, A_2), {\text{gap}} (A_1^{\prime}, A_2^{\prime})\bigr\} \bigm\vert A_1, A_2 \in \alpha,  A^{\prime}_1, A^{\prime}_2 \in \alpha^{\prime}, A_1 \ne A_2 {\text{ or }} A_1^{\prime} \ne A_2^{\prime}\right\}\\
&= \min \left\{\min \bigl\{\max \{{\text{gap}} \left(A_1, A_2\right), {\text{gap}} \left(A_1^{\prime}, A_2^{\prime}\right)\} \bigm\vert A_1 \ne  A_2 \in \alpha, A_1^{\prime} \ne A_2^{\prime} \in \alpha^{\prime}\bigr\},\right.\\
&{\ \  \ }\left.  \min \bigl\{{\text{gap}} \left(A_1, A_2\right) \bigm\vert A_1 \ne A_2 \in \alpha\bigr\}, \min \bigl\{{\text{gap}} \left(A_1^{\prime}, A_2^{\prime}\right) \bigm\vert A_1^{\prime} \ne A_2^{\prime} \in \alpha^{\prime} \bigr\}\right\}\\
&= \min\bigl\{\max \{\gamma (\alpha), \gamma (\alpha^{\prime})\}, \gamma (\alpha), \gamma (\alpha^{\prime})\bigr\}\\
&= \min \bigl\{ \gamma (\alpha), \gamma (\alpha^{\prime})\bigr\}.
\endaligned\]

\noindent Consequently, we have that 
$\max \bigl\{R (\alpha), R (\alpha^{\prime})\bigr\} \le R (\alpha \times \alpha^{\prime})$.

Assume  now  that $\max \bigl\{R (X, d), R (X^{\prime}, d^{\prime})\bigr\} < \infty$.  Let $\bigl(\alpha_n\bigr)_{n \ge1}$ and  $\bigl(\alpha_n^{\prime}\bigr)_{n \ge1}$ be two sequences of clopen partitions of $(X, d)$ and $(X^\prime, d^\prime)$, respectively, with $\delta (\alpha_n) \to 0$,  $\delta (\alpha^{\prime}_n) \to 0$,  $R (\alpha_n) \to R (X, d)$,  and $R (\alpha_n^{\prime}) \to R (X^{\prime}, d^{\prime})$.  Hence, $R (\alpha_n \times \alpha_n^{\prime}) \ge R (X \times X^{\prime}, d_{\sup})$  so that  

\[R (X \times X^{\prime}, d_{\sup})  \le \lim_{k \to \infty} \, \inf_{k \le n} R (\alpha_n \times \alpha^{\prime}_n). \]



\section{Bi-H\"older embedding into Euclidean space}

In this section we consider bi-H\"older embeddings of compact metric spaces of finite logarithmic ratio into $\Bbb R^N$. We are in particular interested in embeddings $f$ with both $f$ and $f^{-1}$ H\"older of same exponent $s$ for some $s > 0$. We also give an estimate for the integer $N$ in terms of the logarithmic ratio and the Assouad dimension $\dim_A$; see \cite{assouad1, assouad2, assouad3} as well as \cite{Luukk, LU-ML}. The following definition for the notion of {\bf metric dimension}  $\dim_m$ is from \cite{mane} and coincides with  Assouad dimension  for totally bounded metric spaces; see also \cite[Definition 3.2 and Remark 3.4]{LU-ML} and \cite{Luukk}.  Let $(X,d)$ be a metric space and, for $0< r_2 <r_1$, let $\mathcal J(r_1, r_2)$ denote the cardinality of a maximal $r_2$-separated set in a closed $r_1$-ball. Then the metric dimension of $X$, denoted by $\dim _m (X, d)$, is defined by setting

\[\dim_m(X,d)\,=\,\lim_{r\to 0}\lim_{t\to \infty} \sup\left\{\frac{\log \mathcal J (r_1, r_2)}{\log (r_1/ r_2)} \biggm\vert 0<r_2 <r_1<r {\text{ and }} r_1 /r_2 >t \right\} .\]

\noindent These limits exist because they are monotone, but they may be infinite.

\vskip5pt

We use this definition because it is more convenient for the proof of Theorem \ref{4.2} below.  It is shown in \cite[Corollary 3.9]{LU-ML} that every ultrametric space of finite Assouad dimension admits a bi-H\"older embedding into $\R^1$. However, that the Assouad dimension be finite is {\it {not}} a necessary condition for a metric space to be bi-H\"older embeddable into $\R^N$ for some $N\in \N$. Indeed, it is shown in \cite[Example 3.6]{LU-ML} that there exists a compact ultrametric space of infinite Assouad dimension which admits a bi-H\"older embedding into $\R^1$. Thus, metric spaces that are bi-H\"older embeddable into $\R^N$ for {\it {some}} $N\in \N$, need {\it {not}} be bi-Lipschitz embeddable into $\R^N$ for {\it {any}} $N\in \N$.

\vskip5pt

Now let $(X,d)$ be a compact metric space of finite logarithmic ratio and assume that $\A^*_p (X) \ne \emptyset$ for some $p \in (1, \infty)$. We know, by Theorem 2.6, that there is a compatible ultrametric $\rho$ on $X$ such that $(X, \rho)$ is bi-H\"older equivalent to $(X, d)$. Then as an immediate corollary  we get the following result.

\begin{cor} Let $(X, d)$ be a compact metric space of finite logarithmic ratio and assume that $\A^*_p (X) \ne \emptyset$ for some $p \in (1, \infty)$. If the ultrametric space $(X, \rho)$ is of finite Assouad dimension, then there is a bi-H\"older embedding of $(X, d)$ into $\R^1$.
\end{cor}

As mentioned above,  the Assouad dimension does not behave well with respect to bi-H\"older maps. Hence, it is possible for $(X, d)$ to be of finite Assouad dimension while $(X, \rho)$ is of infinite Assouad dimension. However, if $\dim_m (X, d) < \infty$ and,    for some $p \in (1, \infty)$, the set $\A^*_p (X)$ contains a very special type of sequence,\footnote{The reason for this extra hypothesis is the  fact that a given proper subsequence of $\bigl(\alpha_n\bigr)_{n\geq 1}\in \A^*_p (X)$ is not necessarily  in $\A^*_q (X)$ for any $q \ge p$.} then there is a bi-H\"older embedding   of $(X, d)$ into $\R^N$ for $N\in \N$ sufficiently large, provided that $1 < R (X, d) < \infty$.

\begin{thm}
\label{4.2}
 Let $(X, d)$ be a compact metric space of Assouad dimension $D < \infty$ and logarithmic ratio $R$ with $1 < R < \infty$.  Suppose that for some $p \in (1, \infty)$ and some $s \in (0, p - 1]$, there is a sequence $\bigl(\alpha_n\bigr)_{n\geq 1}\in \A^*_p (X)$ with $\delta (\alpha_{n + 1}) \le \delta (\alpha_n)^{1 + s}$. Then for any $N > (D + R -1) \bigl[(1 + s) (2R - 1) - 1\bigr]/s$  and any $\varepsilon \in (0, \min \{1, R - 1\})$, there is a bi-H\"older embedding $f: X\hookrightarrow \R^N$  with the property that both $f$ and $f^{-1}$ are H\"older of exponent $1/p(R+\varepsilon)$.
\end{thm}

\noindent In the special case that  $s = R - 1  \in (0, \, p]$, the estimate for $N$ reduces to $N > (D + R -1) (2R+1)$.

\begin{proof}
Let $\bigl(\alpha_n\bigr)_{n\ge 1}\in \A^{\ast}_p (X)$ satisfy the required hypotheses.   We may assume that $\delta (\alpha_n) = \delta_n < 1$ for all $ n \in \N$.  Then, of course, we have that $\delta_n \to 0$ and $\gamma (\alpha_n) = \gamma_n \to 0$ as $n$ tends to infinity with $a \, \delta_n^p \le \delta_{n+1} \le \delta_n^{s + 1}$ for  some $a > 0$ and all $n\ge 1$, and $\lim_{n \to \infty} R (\alpha_n)  =  R (X, d)  =  R$. We note that for each $n \ge 1$ we have 

 \[\max \bigl\{\card (\alpha_{n+1}\cap A)  \bigm\vert \, A \in \alpha_n \bigr\} \le \mathcal J (\delta_n, \gamma_{n+1}),\]

\noindent  where $\alpha_{n + 1} \cap A = \bigl\{A^{\prime} \in \alpha_{n + 1} \bigm\vert A^{\prime} \subset A\bigr\}$, and of course $\delta_n / \gamma_{n+1} \ge \delta_{n+1} / \gamma_{n+1}$ for all $n \geq 1$. Given $\varepsilon > 0$, there exists $n_0$ such that $n_0 \leq n$ implies that $\delta_n^{R+\varepsilon} < \gamma_n < \delta_n^{R-\varepsilon}$. This means that for $n_0 \le n$ we have

\[\delta_{n+1}^{1-(R-\varepsilon)} \, <\, \frac{\delta_{n+1}}{\gamma_{n+1}}\,<\,\delta_{n+1}^{1- (R+\varepsilon)}\,\,.\]

\noindent By choosing  $0 < \ep < R - 1$, we have
 
\[\frac{\delta_n}{\gamma_{n+1}} \,\geq\,\frac{\delta_{n+1}}{\gamma_{n+1}}\,>\,\delta_{n+1}^{1-(R-\varepsilon)}\]

\noindent which approaches infinity as $n$ tends to infinity. Thus, by definition, there exist constants $n_1 \ge n_0$ and $K (\varepsilon)$ such that $n_1 \le n$ and $K (\varepsilon) < \delta_n /\gamma_{n+1}$ imply that

\[\mathcal J (\delta_n, \gamma_{n+1})\,\leq\,\left(\frac{\delta_n}{\gamma_{n+1}}\right)^{D+\varepsilon} .\]

\noindent Therefore, if $n_1 \le n$ and $K (\varepsilon) < \delta_n /\gamma_{n+1}$, then 

\begin{equation}
\label{4.2.1}
\max \bigl\{ \card (\alpha_{n+1}\cap A) \, \bigm\vert \, A\in\alpha_n\bigr\} \, \leq \, \left(\frac{\delta_n}{\gamma_{n+1}}\right)^{D+\varepsilon} .
\end{equation}

The rest of the proof consists of two steps.

\vskip5pt

\noindent {\it {Step 1.}} By using the box norm $\Vert \cdot\Vert$ in $\R^N$, it is easily seen that in a closed ball of radius $\delta_n$  we can fit 

\[\left\lfloor \frac{\delta_n + \gamma_{n +1}/2}{\delta_{n+1} + \gamma_{n+1}/2}\right\rfloor^N\]

\noindent disjoint closed balls of radius $\delta_{n+1}$ and gap $\gamma_{n+1}$. Thus, we need to have

\[\max \bigl\{ \card (\alpha_{n+1}\cap A)\, \bigm\vert \, A\in\alpha_n\bigr\}\,\leq\,\left\lfloor\frac{\delta_n + \gamma_{n +1}/2}{\delta_{n+1} + \gamma_{n+1}/2}\right\rfloor^N .\]

\noindent By (\ref{4.2.1}), it suffices to have

\begin{equation}
\label{4.2.2}
\left(\frac{\delta_n}{\gamma_{n+1}}\right)^{D+\varepsilon}\,\leq\,\left\lfloor\frac{\delta_n + \gamma_{n +1}/2}{\delta_{n+1} + \gamma_{n+1}/2}\right\rfloor^N .
\end{equation}

\noindent Inequality (\ref{4.2.2}) holds provided that for $n$ sufficiently large we have

\[\delta_{n+1}\,<\,\delta_n^{[N-(D + \ep)] / [N - (R +\ep) (D + \ep)]}.\]

\noindent But, since  $0 < \ep < R -1$ and $N > (D + R -1) \bigl[(1 + s) (2R - 1) - 1\bigr]/s$, we see that 

\[1 + s > \frac{N - (D + \ep)}{N - (R + \ep)(D + \ep)}\]

\noindent  so that  for $n$ sufficiently large we have

\[a \,  \delta_n^p \le \delta_{n+1} \le  \delta_n^{1 + s} \, < \, \delta_n^{[N-(D + \ep)]/[N-(R + \ep)(D + \ep)]}\]

\noindent as needed. 

\vskip5pt

\noindent {\it Step 2.} The following construction is a variant of the one in \cite{mane}. For each $n\geq 1$, define  $g_n : \alpha_{n+1}\rightarrow\alpha_n$ by setting $g_n (A)\,=$ the member of $\alpha_n$ containing $A;\, A\in \alpha_{n+1}$. In $\R^N$ let $\mathcal B_1$ denote a family of $(\card\,\alpha_1)$-many disjoint closed balls of radius $\delta_1$ with $\text{dist}(B^{\prime}, B^{\prime\prime}) \geq \gamma_1$ for $B^{\prime} \neq B^{\prime\prime}$. By (\ref{4.2.2}), inside each $B\in\mathcal B_1$ we may choose $\card (\alpha_2 \cap B)$-many disjoint closed balls of radius $\delta_2$ with $\text{dist}(B^{\prime}, B^{\prime\prime}) \geq \gamma_2$ for  $B^{\prime} \neq B^{\prime\prime}$. Let $\mathcal B_2$ denote the family of these balls. By induction, if the family $\mathcal B_n$ has been chosen, inside each ball $B\in \mathcal B_n$ we may choose $\card (\alpha_{n+1} \cap B)$-many disjoint closed balls of radius $\delta_{n+1}$ with the property that $\text{dist}(B',B'')\geq \gamma_{n+1}$ for $B'\neq B''$.
We let $\mathcal B_{n+1}$ denote the family of these balls and the induction is complete.

\vskip5pt

 Next, for each $n\geq 1$ let $h_n : \mathcal B_{n+1}\rightarrow \mathcal B_n$ by setting $h_n (B)\,=$ the member of $\mathcal B_n$ containing $B;\, B\in\mathcal B_{n+1}$. Let $f_1 : \alpha_1\rightarrow\mathcal B_1$ be any bijection and choose a bijection $f_2 : \alpha_2\rightarrow\mathcal B_2$ so that $f_1 g_1\,=\,h_1 f_2$. Similarly, if $f_n :\alpha_n\rightarrow\mathcal B_n$ has been chosen, we may choose a bijection $f_{n+1} : \alpha_{n+1}\rightarrow\mathcal B_{n+1}$ so that $f_n g_n\,=\,h_n f_{n+1}$. In short, we have the following commutative diagram:
        
\[\xymatrix{
\alpha_1 \ar[d]_{f_1} &\alpha_2 \ar[l]_{g_1} \ar[d]_{f_2} &\alpha_3 \ar[l]_{g_2} \ar[d]_{f_3} &\cdots \ar[l]_{g_3}\\
\mathcal B_1 &\mathcal B_2 \ar[l]^{h_1} &\mathcal B_3 \ar[l]^{h_2} &\cdots \ar[l]^{h_3}
}\]

\noindent Finally, define $f : X \aro \R^N$ by setting

\[\{f(x)\} \,=\,\bigcap_{n\geq1} f_n\left(\alpha_n(x)\right) ,\]

\noindent where $\alpha_n(x)\,=$ the member of $\alpha_n$ containing $x$. If $x_1,x_2\,\in X$ are in the same member of $\alpha_n$ but in different members of $\alpha_{n+1}$, then $\gamma_{n+1} \leq d(x_1,x_2) \leq \delta_n$. Furthermore, $f(x_1)$ and $f(x_2)$ are in the same member of $\mathcal B_n$ but in different members of $\mathcal B_{n+1}$ implying that $\gamma_{n+1} \leq\,\,\parallel f(x_1)-f(x_2) \parallel \,\,\leq 2\,\delta_n$. But, there exists $a > 0$ with

\[a \, \delta_n^p\,<\,\delta_{n+1}\,<\,\gamma_{n+1}^{1/(R+\varepsilon)}\,\leq\,d(x_1,x_2)^{1/(R+\varepsilon)}\]

\noindent so that $\delta_n <  a^{- 1/p} d(x_1, x_2)^{1/p(R+\varepsilon)}$. Moreover, we have $\gamma_{n+1} > \delta_{n+1}^{R+\varepsilon} \geq  a^{R + \varepsilon} \delta_n^{p(R+\varepsilon)}$ which gives $\gamma_{n+1} \geq d(x_1, x_2)^{p(R+\varepsilon)}$. Hence,

\[a^{R + \varepsilon} d(x_1,x_2 )^{p(R+\varepsilon)} \leq \bigl\Vert f(x_1) - f(x_2)\bigr\Vert \leq 2 a^{- 1/p}\, d(x_1, x_2)^{1/p(R+\varepsilon)}.\] 
\end{proof}



\section{Remarks}

\noindent {\bf{I.}} For comparison with Corollary 2.2, we note  that by \cite[Theorem A.5]{Luukk} we have $\dim_A (X,d^s) = (1/s) \dim_A (X, d)$.

\vskip5pt
 
\noindent {\bf{II.}} If Conjecture 2.8 is true, then Lemma 2.5 and Theorem 4.2 are strengthened significantly. In particular, Theorem \ref{ult} would then imply that the logarithmic ratio  gives a complete characterization of {\it all} compact metric spaces which are bi-H\"older equivalent to ultrametric spaces.

\vskip5pt

\noindent{\bf{III.}}   There exists a compact ultrametric space $(X, d)$ homeomorphic to the Cantor ternary subset of $[0, 1]$ for which $R (X, d) = 0$.    

\begin{examples}   Let $X = \left(\Z /2\Z\right)^{\N}$ and equipped  with  the ultrametric $d$ defined as follows: Let $0 < r < 1$ be fixed and let $r_n = r^{(n!)}$. For $x = (x_1, x_2, \cdots)$ and $x^{\prime} = (x_1^{\prime}, x_2^{\prime}, \cdots)$ in $X$, let $d(x, x^{\prime}) = r_n$ if $x_i = x_i^{\prime}$ for $1 \le i \le n$ and $x_{n + 1} \ne x_{n + 1}^{\prime}$. Now for each $n > 1$, let $\alpha_n$ be the partition of $X$ into $r_n$--balls. Then $\delta (\alpha_n) = r_n$ and $\gamma (\alpha_n) = r_{n - 1}$ and hence $R(\alpha_n) = 1/n$. Therefore, $0 \le R(X, d) \le \liminf_{n \to \infty} R (\alpha_n) = 0$. We note that as in  Example 3.1, Case 3.1.1,  we have  $\delta (\alpha_{n + 1})/ \delta (\alpha_n)^p \geq 1$ if and only if $p \geq n + 1$. Hence, for each $p > 1$ and each $n(p) \in \N$ the sequence $\{\alpha_n\}_{n\geq n(p)} \notin {\A}^*_p(X)$.   $\square$
\end{examples}

\vskip5pt

\noindent {\bf{IV.}} Although the case $R(X,d) = 1$ appears to be anomalous (see, for instance, \cite[Propositions 2.12, 2.13]{lograt}   and \cite{ultgeomeas}), it is still possible to obtain partial results (in spirit of Theorem 4.2) in this case.

\vskip5pt

\noindent {\bf{V.}} It is noted in \cite{lograt} that, in general, it is not easy to compute the logarithmic ratio of a compact (totally disconnected, nondiscrete) metric space $(X, d)$. Indeed,  most of our computations use \cite[Theorem 2.7]{lograt} since it is the only computational result we have.  However, we obtain upper and lower bounds for the logarithmic ratio as follows:  We define two functions $ G, g : [0, \diam X] \aro [0, \diam X]$ by setting 

\[\aligned
G (r) &= \inf\{\gamma (\alpha) \bigm\vert \alpha {\text{ is a clopen partition of }} X {\text{ with }} \delta (\alpha) \ge r\}{\text{ and }}\\
g (r) &= \sup\{\gamma (\alpha) \bigm\vert \alpha {\text{ is a clopen partition of }} X {\text{ with }} \delta (\alpha) \le r\}.
\endaligned\] 

\noindent Clearly, both functions $G$ and $g$ are monotone increasing with        $\lim_{r \to 0} G (r) = 0$ and  $\lim_{r \to 0} g (r) = 0$.   Let

\[\underline{\mathcal R} (X, d) = \liminf_{r \to 0} \frac{\log g (r)}{\log r} \qquad {\text{ and }} \qquad \overline{\mathcal  R} (X, d) = \liminf_{r \to 0} \frac{\log G (r)}{\log r}.\]

\noindent Then we have

\[
\underline{\mathcal R} (X, d) \le R (X, d) \le \overline{\mathcal R} (X, d).\]

\vskip5pt

\noindent {\bf VI. } In Theorem 4.2, we have $R / p^2 (R + \ep)^2 \le R \bigl(f (X), \Vert . \Vert\bigr) \le R  p^2 (R + \ep)^2$ and $\dim_H (X) /p (R + \ep) \le \dim_H \bigl(f (X)\bigr) \le p (R + \ep) \dim_H (X)$.

\vskip5pt

\noindent {\bf{VII.}}  Let $(X, d)$ be a compact metric space and let $f : X \aro \{0\}$.  Then, obviously, $f$ is Lipschitz with $0 = R (\{0\}) \le R (X, d)$.  However, consider  the following question:   Can a Lipschitz map of compact metric spaces raise the logarithmic ratio?   The following example, due to an anonymous reader, shows that the answer is affirmative. 

\begin{lemma} {\rm{[anonymous reader]}} 
There exists a Lipschitz surjection $f : (X, d) \aro (X^{\prime}, d^{\prime})$ of compact metric spaces with $R (X, d) < R (X^{\prime}, d^{\prime})$.
\end{lemma}

\begin{proof}  Let $X = \{0\} \cup \{1/n \bigm\vert n \in \N\}$ with the topologically compatible ultrametric $d (x, y) = \max \{\sqrt x, \sqrt y\}$.  For each $n > 1$, let $A_n = X \setminus \{1/k \bigm\vert 1 \le k \le n -1\}$ and, as in Example 3.1, let  the partition $\alpha_n = \bigl\{A_n, \{1\}, \{1/2\}, \cdots , \{1/ (n - 1)\}\bigr\}$.   Then $\delta (\alpha_n) = 1/\sqrt n$ and  $\gamma (\alpha_n) = 1 / \sqrt{n - 1}$. Also $\diam (A_n) = \Gamma (A_n) = 1 / \sqrt n = \gamma (\alpha_{n + 1})$.  Hence, 
\[R (X, d) = \liminf_{n \to \infty} R (\alpha_n) = \lim_{n \to \infty} \frac{\log (n - 1)}{\log n} = 1.\]

\noindent Moreover, if $x, y \in X$ with $x > y$, then $\vert x - y\vert = x - y \le x \le \sqrt x = d (x, y)$ so that the identity map $\id : (X, d) \aro (X, \vert \cdot \vert)$ is $1$-Lipschitz.  However, by the case $s = 2$ of Case 3.2.4 in  Example 3.1 we have that $R (X, \vert \cdot \vert)= 2$.
 \end{proof}

\vskip5pt

\noindent {\bf VIII.}    For a metric space $(X, d)$, let $\h(X)$ denote the space of all compact nonempty subsets of $(X, d)$ equipped with the Hausdorff metric. The topological structure of $\h  (X)$ has been studied extensively; see \cite{Illanes, schoriwest} and references therein.  For instance, it is known that if $(X, d)$ is complete (respectively compact), then $\h (X)$ is also complete (respectively compact). Moreover, if $(X, d)$ is separable, then so is $\h (X)$.   Indeed, it is a celebrated result of  Schori and  West \cite{schoriwest} that $\h ([0, 1])$ is homeomorphic to the Hilbert cube $[0, 1]^{\N}$.

\vskip5pt

The metric structure of $\h (X)$ has also been studied extensively \cite{bandt, boardman, goodey1, goodey2, Grub, Grub-Let, Grub-Tich, Hohti, McClure, tyson}.  For example it is shown in \cite{Hohti} that $\h ([0, 1])$ is not bi-Lipschitz equivalent to a variety of metric Hilbert cubes.  The intrinsic metric complexity of $\h ([0, 1])$ is further illustrated in \cite{tyson} by showing that there is {\it{no}} bi-Lipschitz embedding of $\h ([0, 1])$ into any uniformly convex Banach space.  However, there are some positive results. Specifically, it is immediate from the definition that if $(X, d)$ is an ultrametric space, then $\h (X)$ is also an ultrametric space.  Since every separable  ultrametric space admits an isometric embedding into a Hilbert space \cite{aschbacher, kelley, Lemin, timan, timvest}  and since a bi-H\"older equivalence of compact metric spaces induces a bi-H\"older equivalence of their hyperspaces, the following result is an immediate corollary of  Theorem 2.6.

\begin{cor} Let $(X, d)$ be a compact metric space of finite logarithmic ratio and assume that ${\A}^*_p (X)\ne\emptyset$ for some $p\in (1, \infty)$. Then there is a bi-H\"older embedding of  $\h (X)$ into a Hilbert space.
\end{cor}

\begin{examples}  Let $(X, d)$ be  the metric space in Example 3.1, Case 3.1.1.  It   follows from     \cite[Theorem 2.8]{tyson} along with straightforward calculations  that  there is a bi-H\"older (in fact, bi-Lipschitz) embedding $f : \h (X) \hookrightarrow \R^2$.    
\end{examples}



\end{document}